\documentclass[11pt,twoside,a4paper,leqno]{amsart}
\calclayout
\usepackage{amssymb}
\newif\ifaddpics\addpicstrue
\ifaddpics\usepackage{graphicx}\fi
\makeatletter
\def\swappedhead@plain#1#2#3{%
  \thmnumber{\@upn{\mdseries #2}}\thmname{\@ifnotempty{#2}{. }#1}%
  \thmnote{ \textmd{\upshape(#3)}}}
\makeatother
\theoremstyle{plain}
\newtheorem{result}{Theorem}

\swapnumbers
\newtheorem{thm}[subsection]{Theorem}
\newtheorem{prop}[subsection]{Proposition}
\theoremstyle{remark}
\newtheorem*{rems}{Remarks}
\newtheorem*{examples}{Examples}

\newcommand{\emphdef}{\textit}                  
\newcounter{numl}

\newcommand{\labelnuml}{\textup{(\roman{numl})}}
\newenvironment{numlist}{\begin{list}{\labelnuml}%
{\usecounter{numl}\setlength{\leftmargin}{0pt}%
\setlength{\itemindent}{2\parindent}%
\setlength{\itemsep}{\smallskipamount}\def
\makelabel ##1{\hss \llap {\upshape ##1}}}}{\end{list}}
\newenvironment{bulletlist}{\begin{list}{\labelitemi}%
{\setlength{\leftmargin}{\parindent}\def
\makelabel ##1{\hss \llap {\upshape ##1}}}}{\end{list}}
\newcommand{\acknowledge}{\subsection*{Acknowledgements}}
\newcommand{\thismonth}{\ifcase\month\or
  January\or February\or March\or April\or May\or June\or
  July\or August\or September\or October\or November\or December\fi
  \space\number\year}
\makeatletter
\newcommand{\low}{\@ifnextchar^{}{^{\vphantom x}}}
\newcommand{\high}{\@ifnextchar_{}{_{\vphantom I}}}
\makeatother
\DeclareSymbolFont{script}{U}{eus}{m}{n}
\DeclareSymbolFontAlphabet{\mathscr}{script}
\DeclareMathSymbol{\EuWedge}{0}{script}{"5E}
\DeclareMathAlphabet{\mathrmsl}{OT1}{cmr}{m}{sl}
\newcommand{\rssymb}[2]{\newcommand{#1}{{\mathrmsl{#2}}}}
\newcommand{\calsymb}[2]{\newcommand{#1}{{\mathcal{#2}}}}
\newcommand{\bbsymb}[2]{\newcommand{#1}{{\mathbb{#2}}}}
\newcommand{\lieoper}[2]{\newcommand{#1}{\mathop
  {\mathfrak{#2}\null}\nolimits}}
\newcommand{\oper}[3][n]{\newcommand{#2}{\mathop
  {\mathrm{#3}\null}\ifx n#1\nolimits\else\limits\fi}}
\newcommand{\rsoper}[3][n]{\newcommand{#2}{\mathop
  {\mathrmsl{#3}\null}\ifx n#1\nolimits\else\limits\fi}}
\bbsymb\C{C} \bbsymb\F{F} \bbsymb\HQ{H}\bbsymb\N{N} \bbsymb\IP{P}
\bbsymb\Q{Q} \bbsymb\R{R} \bbsymb\U{U} \bbsymb\V{V} \bbsymb\W{W} \bbsymb\Z{Z}
\calsymb\cA{A} \calsymb\cB{B} \calsymb\cC{C} \calsymb\cD{D} \calsymb\cE{E}
\calsymb\cF{F} \calsymb\cG{G} \calsymb\cH{H} \calsymb\cI{I} \calsymb\cJ{J}
\calsymb\cK{K} \calsymb\cL{L} \calsymb\cM{M} \calsymb\cN{N} \calsymb\cO{O}
\calsymb\cP{P} \calsymb\cQ{Q} \calsymb\cR{R} \calsymb\cS{S} \calsymb\cT{T}
\calsymb\cU{U} \calsymb\cV{V} \calsymb\cW{W} \calsymb\cX{X} \calsymb\cY{Y}
\calsymb\cZ{Z}
\newcommand{\eps}{\varepsilon}

\renewcommand{\geq}{\geqslant} \renewcommand{\leq}{\leqslant}
\oper\End{End} \oper\Hom{Hom}                    
\oper\Sym{Sym} \oper\Skew{Skew}
\oper\Aut{Aut}                                   
\oper\GL{GL} \oper\SL{SL}\oper\Symp{Sp}
\oper\CO{CO} \oper\On{O} \oper\SO{SO} \oper\Pin{Pin} \oper\Spin{Spin}
\oper\CU{CU} \oper\Un{U} \oper\SU{SU}
\rsoper\Diff{Diff} \rsoper\SDiff{SDiff}
\lieoper\der{der}                                
\lieoper\gl{gl} \lieoper\sgl{sl}\lieoper\symp{sp}
\lieoper\co{co} \lieoper\so{so} \lieoper\spin{spin}
\lieoper\cu{cu} \lieoper\un{u}  \lieoper\su{su}
\rsoper\Vect{Vect} \rsoper\Ham{Ham}
\newcommand{\norm}[2][]{|\mkern-2mu|#2|\mkern-2mu|
  _{\lower1pt\hbox{${}_{#1}$}}}
\newcommand{\Norm}[2][]{\bigl|\mkern-3mu\bigr|#2\bigr|\mkern-3mu\bigr|
  _{\lower1pt\hbox{${}_{#1}$}}}

\newcommand{\punc}[1]{\smallsetminus\{#1\}}

\newcommand{\tens}{\otimes}                 
\newcommand{\del}{\partial}                 
\newcommand{\Proj}{\mathrmsl{P}}            
\newcommand{\RP}[1]{\R\Proj^{#1}}           
\newcommand{\HP}[1]{\HQ\Proj^{#1}}          
\rsoper\dimn{dim}                           
\rsoper\grad{grad}                          
\rsoper\kernel{ker}\rsoper\image{im}        
\rsoper\alt{alt}   \rsoper\sym{sym}         
\rsoper\Ad{Ad}     \rsoper\ad{ad}           
\rsoper\CoAd{CoAd} \rsoper\coad{coad}       
\rsoper\trace{tr}  \rsoper\trfree{tf}       
\rsoper\detm{det}                           
\rsoper\Vol{Vol}                            
\rsoper\divg{div}                           
\rsoper\sign{sign}
\rssymb\iden{id}                            
\rssymb\vol{vol}                            
\oper\Imag{Im}
\newcommand{\ips}[1]{\eps(#1)}

\newcommand{\sd}{{\raise1pt\hbox{$\scriptscriptstyle +$}}}
\newcommand{\asd}{{\raise1pt\hbox{$\scriptscriptstyle -$}}}
\newcommand{\sdasd}{{\raise1pt\hbox{$\scriptscriptstyle\pm$}}}
\newcommand{\asdsd}{{\raise1pt\hbox{$\scriptscriptstyle\mp$}}}

\rsoper\scal{scal}

\begin{document}
\title{Toric selfdual Einstein metrics on compact orbifolds}
\author{David M. J. Calderbank}
\address{School of Mathematics\\
University of Edinburgh\\ King's Buildings, Mayfield Road\\ Edinburgh
EH9 3JZ\\ Scotland.}
\email{davidmjc@maths.ed.ac.uk}
\author{Michael A. Singer}
\email{michael@maths.ed.ac.uk}
\date{\thismonth}
\begin{abstract}
We prove that any compact selfdual Einstein $4$-orbifold of positive scalar
curvature whose isometry group contains a $2$-torus is, up to an orbifold
covering, a quaternion K\"ahler quotient of $(k-1)$-dimensional quaternionic
projective space by a $(k-2)$-torus for some $k\geq 2$.  We also obtain a
topological classification in terms of the intersection form of the
$4$-orbifold.
\end{abstract}

\maketitle
\section*{Introduction}

A selfdual Einstein (SDE) metric is a $4$-dimensional Riemannian metric $g$
whose Weyl curvature $W$ is selfdual with respect to the Hodge star
operator ($W=*W$), and whose Ricci tensor is proportional to the metric
(Ric $=\lambda g$). The only compact oriented $4$-manifolds admitting SDE
metrics of positive scalar curvature are $S^4$ and $\C P^2$, with the round
metric and Fubini--Study metric respectively (cf.~\cite[Thm.~13.30]{Bes:em}).
However, if one considers $4$-orbifolds, the class of examples is much
wider. In~\cite{GaLa:qro}, K. Galicki and H.~B. Lawson constructed SDE
$4$-orbifolds by taking quaternion K\"ahler quotients of quaternionic
projective spaces by tori, and this construction was systematically
investigated by Boyer--Galicki--Mann--Rees~\cite{BGMR:3s7}.

These examples are all \emph{toric}, i.e., the isometry group of the metric
contains a $2$-torus, hence they belong to the local classification by
H. Pedersen and the first author of toric SDE metrics of nonzero scalar
curvature~\cite{CaPe:emt}, where it was shown that any such metric has an
explicit local form determined by an eigenfunction of the Laplacian on the
hyperbolic plane $\cH^2$.

SDE $4$-orbifolds have physical relevance as the simplest nontrivial target
spaces for nonlinear sigma models in $N=2$ supergravity
(see~\cite{Gal:nmc}). They have also attracted interest recently because of
the connection with M-theory and manifolds of
$G_2$-holonomy~\cite{BrSa:ceh,GPP:emb}: in particular the results in
\cite{CaPe:emt} have been exploited by L. Anguelova and
C. Lazaroiu~\cite{AnLa:mtg2}.

The first main theorem of this paper shows that the quaterion K\"ahler
quotient is sufficient.

\begin{result}\label{thma}
Let $X$ be a compact selfdual Einstein $4$-orbifold with positive scalar
curvature, whose isometry group contains a $2$-torus. Then, up to orbifold
coverings, $X$ is isometric to a quaternion K\"ahler quotient of
quaternionic projective space $\HP{k-1}$, for some $k\geq 2$, by a
$(k-2)$-dimensional subtorus of $\Symp(k)$. \textup(We remark that the
least such $k$ is $b_2(X)+2$.\textup)
\end{result}

This result was already known when the $3$-Sasakian $7$-orbifold associated
to $X$ is smooth, since toric $3$-Sasakian $7$-manifolds have been
classified by R. Bielawski~\cite{Bie:chkt} using analytical techniques. Our
methods are quite different, being elementary and entirely explicit.

Before outlining the proof, we recall that it was shown in \cite{CaPe:emt}
that the SDE metrics coming from (local) quaternion K\"ahler quotients of
$\HP{k-1}$ as above, are those for which the corresponding hyperbolic
eigenfunction $F$ is a positive linear superposition of $k$ basic solutions
which may be written
\begin{equation}\label{multipole}
F(\rho,\eta)= \sum_{i=1}^k
\frac{\sqrt{a_{\smash i}^{2}\rho^2+(a_i\eta-b_i)^2}}{\sqrt{\rho}},
\end{equation}
where $a_i,b_i \in \R$, and $(\rho>0,\eta)$ are half-space coordinates on
$\cH^2$. Our main task, therefore, is to show that the eigenfunction $F$
associated to $X$ in Theorem~\ref{thma} is of the form \eqref{multipole}.

\begin{proof} We provide here the main line of the argument,
relegating the detail to the body of the paper. First we observe that by
Myers' Theorem (which extends easily to orbifolds~\cite{Bor:omd}) the
universal orbifold cover $\tilde X$ of $X$ is also a compact toric SDE
$4$-orbifold of positive scalar curvature, so we may assume $X=\tilde X$.

As we shall explain in section 1, compact simply connected $4$-orbifolds $X$
with an effective action of a $2$-torus $G=T^2$ may be classified by work of
Orlik--Raymond~\cite{OrRa:at4} and Haefliger--Salem~\cite{HaSa:ato}.  The
orbit space $W=X/G$ is a polygonal disc whose edges $C_0,C_1,\ldots
C_{k-1},C_k=C_0$ (given in cyclic order) are labelled by \emphdef{orbifold
generators} $v_j=(m_j,n_j)\in \Z^2$, determined up to sign\footnote{Thus
$v_k=\pm v_0$ and it will be convenient later to take the sign to be
negative.}, with $m_j n_{j-1}-m_{j-1}n_j\neq 0$ for $j=1,\ldots k$. The
interior of $W$ is the image of the open subset $X_0$ of $X$ on which $G$
acts freely, the edges $C_j$ are the images of points with stabilizer
$G(v_j)=\{(z_1,z_2)\in G: m_j z_1 + n_j z_2=0\}$ and cyclic orbifold
structure groups of order $\mathrm{gcd}(m_j,n_j)$, and the corners are the
images of the fixed points.  A sign choice for $v_j$ is equivalent to an
orientation of the corresponding circle orbits.

The classification result of~\cite{CaPe:emt}, which we discuss in section
2, shows that the interior of $W$ is equipped with a hyperbolic metric
$g_{\cH^2}$ and hyperbolic eigenfunction $F$ (with $\Delta_{\cH^2} F=\frac
34 F$) such that the SDE metric $g_F$ on $X_0$ is given explicitly by the
formula~\eqref{metric}.

We shall show in section 3 that for the compactification of the metric
$g_F$ on $X$, it is necessary that the edges of $W$ are at infinity with
respect to the hyperbolic metric. It follows (by simple connectivity) that
the interior of $W$ is identified with the entire hyperbolic plane.  We
also show in section 3 that for any half-space coordinates $(\rho>0,\eta)$
on $\cH^2$, the function $\sqrt{\rho}F(\rho,\eta)$ has a well-defined limit
as $\rho\to 0$, which is a continuous piecewise linear function $f_0(\eta)$
of $\eta$ (whose corners are at the corners of $W$). The half-space
coordinates can be chosen so that $f_0(\eta)= \pm(m_j\eta-n_j)$ on $C_j$.

As we shall discuss in section 2 (cf.~\cite{CaSi:emcs}) a hyperbolic
monopole $F$ on $\cH^2$ is determined by its `boundary value' $f_0$ via a
`Poisson formula'
\begin{equation*}
F(\rho,\eta)=\frac{1}{2}\int \frac{f_0(y)\rho^{3/2}\,d y}{\bigl(\rho^2 + (\eta
-y)^2\bigr)^{3/2}}.
\end{equation*}
Integrating twice by parts (in the sense of distributions), we then have
\begin{equation*}
F(\rho,\eta)=\frac{1}{2}\int \frac{f_0''(y)
\sqrt{\rho^2 + (\eta -y)^2}\,d y} {\sqrt{\rho}}.
\end{equation*}
In our case, $f_0$ is continuous and piecewise linear, so $f_0''$ is a
linear combination of $k$ delta distributions and $F$ is therefore a linear
combination of $k$ basic solutions, i.e., a $k$-pole solution in the sense
of~\cite{CaPe:emt}. Since the SDE metric has positive scalar curvature, it
follows from~\cite{CaPe:emt} that the determinant of a certain matrix
$\Phi(\rho,\eta)$ associated to $F(\rho,\eta)$ (see section 2) is positive.
Using~\cite{CaSi:emcs} (see section 2 again) we find that
\begin{equation}\label{lovelyphieqn}
\det\Phi(\rho,\eta) = \frac{1}{4}\iint
\frac{f_0(y)f_0''(z)\rho\bigl(\rho^2+(\eta-y)(\eta-z)\bigr)}
{\bigl(\rho^2+(\eta-y)^2\bigr){}^{3/2}\bigl(\rho^2+(\eta-z)^2\bigr){}^{3/2}}
\,dy\,dz.
\end{equation}
Suppose now $\eta$ lies in the singular set of $f_0''$. Then as $\rho\to 0$,
this integral is dominated by the contribution from $y$ near $\eta$ and the
evaluation at $z=\eta$ given by the corresponding delta distribution in
$f_0''$. It follows that $f_0$ is convex where it is positive and concave
where it is negative. Thus, up to an irrelevant sign, $f_0(\eta)$ is
positive and convex, hence $F(\rho,\eta)$ is of the form~\eqref{multipole}.
According to~\cite{CaPe:emt}, the metric $g$ is therefore locally isometric
to a local quaternion K\"ahler quotient of $\HP{k-1}$ by an explicitly
defined $(k-2)$-dimensional abelian subgroup of $\Symp(k)$, and one easily
sees that the integrality conditions $v_j\in \Z^2$ imply that this subgroup
is a torus. However, the quaternion K\"ahler quotient of $\HP{k-1}$ by a
$(k-2)$-torus is a compact $4$-orbifold~\cite{BGMR:3s7}.  Therefore $X$ must
be its universal orbifold cover.
\end{proof}

The proof of this theorem shows explicitly how the isotropy data of a toric
SDE orbifold $X$ give rise to the hyperbolic eigenfunction defining the SDE
metric on $X$ and hence to its realization as quaternion K\"ahler quotient.
However, it is not yet clear which isotropy data give a toric orbifold
admitting an SDE metric. To understand this, we consider the inverse
construction. Suppose therefore that $X$ is a quaternion K\"ahler quotient
of $\HP{k-1}$ by a $(k-2)$-dimensional subtorus $H$ of the standard maximal
torus $T^{k}=\R^k/2\pi \Z^k$ in $\Symp(k)$. Then the isometry group of $X$
contains the quotient torus $T^{k}/H$. If we choose an identification of
$T^k/H$ with with $\R^2/2\pi \Z^2$, then $H=\mathfrak h/2\pi\Lambda$ is
determined by a map from $\Z^{k}\to \Z^2$ with kernel $\Lambda$, or
equivalently by $(a_i,b_i)\in \Z^2$ for $i=1,\ldots k$ (the images of the
standard basis elements of $\Z^k$). Any two $(a_i,b_i)$ span
$\Z^2\tens_{\Z}\Q$: otherwise $X$ is a quaternion K\"ahler quotient of
$\HP{j-1}$ for some $j<k$.  Using the choice of basis of $\Z^2$, we can
suppose $a_i\neq 0$ for all $i$, and then using the ordering and signs of
the standard basis of $\Z^k$, we can assume that $a_i>0$ and that the
sequence $(y_i:=b_i/a_i)$ is increasing. We set $y_0=-\infty$ and
$y_{k+1}=+\infty$.

Now, by~\cite{CaPe:emt}, the hyperbolic eigenfunction generating the SDE
metric is given by~\eqref{multipole} and therefore the boundary value of
$\sqrt{\rho} F(\rho,\eta)$ is the continuous piecewise linear convex
function
\begin{equation}\label{multibdy}
f_0(\eta)={\textstyle\sum}_{i=1}^k |a_i\eta-b_i|,
\end{equation}
whose value on the edge $(y_i,y_{i+1})$ is $m_i\eta-n_i$, where
\begin{align*}
v_j=(m_j,n_j)&={\textstyle\sum}_{i=1}^j (a_i,b_i)
-{\textstyle\sum}_{i=j+1}^k (a_i,b_i)\\
2(a_i,b_i) &= (m_i,n_i) - (m_{i-1},n_{i-1})
\end{align*}
and we set $(m_0,n_0)=-(m_k,n_k)$.  Up to orbifold covering, $X$ is the
compact toric orbifold with stabilizers and orbifold structure groups
determined by $\pm v_j$ in that cyclic order.

In terms of the $v_j$, the conditions $a_i>0$ and $(y_i)$ strictly
increasing mean that:
\begin{numlist}
\item[(a)] the sequence $m_j$ is strictly increasing\textup;
\item[(b)] the sequence $(n_{j}-n_{j-1})/(m_{j}-m_{j-1})$ is strictly
increasing.
\end{numlist}
(Relative to a given basis of $\Z^2$, condition (a) determines the cyclic
permutation of the $v_j$, since it forces $m_0=-m_k<m_j<m_k$ for $0<j<k$,
so $|m_0|=|m_k|$ is the largest $|m_j|$. It also fixes the signs, apart
from that of the smallest $|m_j|$, but this is fixed by condition (b).)

Conversely, given isotropy data for a toric orbifold $X$, if we can choose
the signs and the cyclic permutation so that (a)--(b) hold, we can define
the corresponding $(a_i,b_i)$ and $f_0(\eta)$ inducing these data (up to a
factor of $2$). $X$ will then admit a toric SDE metric since, up to orbifold
covering, it is a quaternion K\"ahler quotient as above.

We now give a topological interpretation of these conditions, which leads to
a classification result (cf.~\cite{AnLa:mtc}). Before stating it, we define
$\Delta_{i,j} = m_in_j - m_jn_i$ and introduce the term \emphdef{exceptional
surface} to refer to the inverse image in a toric orbifold $X$ of an edge of
$X/G$.

\begin{result} Let $X$ be a compact, simply connected, oriented toric
$4$-orbifold, with $k=b_2(X)+2$. Then the following are equivalent.
\begin{numlist}

\item $X$ admits a selfdual Einstein metric of positive scalar curvature
such that the exceptional surfaces are totally geodesic.

\item The intersection form of $X$ is positive definite and for any
exceptional surface $S$, the rational homology class $[\bar S]$ has
self-intersection number $e=[\bar S]\cdot[\bar S]<\chi_{\mathrm{orb}} (\bar
S)$, where $\chi_{\mathrm{orb}} (\bar S)$ is the orbifold Euler
characteristic of the closure $\bar S$ of $S$.

\item If $S_1,S_2,\ldots S_k$ are the exceptional surfaces oriented so that
the orbifold generators $v_1,v_2\ldots v_k$ satisfy $\Delta_{0,j}\geq 0$
where $v_0=-v_k$, then the $v_j$ are in cyclic order
\textup($\Delta_{j-1,j}>0$\textup) and
\begin{equation}\label{isn-formula}
\Delta_{j-1,j+1} < \Delta_{j-1,j} + \Delta_{j,j+1} \qquad \text{for all}
\quad 0<j<k.
\end{equation}

\end{numlist}

If \textup{(i)--(iii)} hold, then $X$ admits a toric selfdual Einstein
metric, unique up to homothety and pullback by an equivariant
diffeomorphism.
\end{result}
\begin{proof} 
(i)$\Rightarrow$(ii).  Since the metric is selfdual with positive scalar
curvature, the intersection form must be positive definite (by Hodge theory
for orbifolds and a Bochner argument---cf.~\cite{LeBr:tsd}). Now we
observe that $e<\chi_{\mathrm{orb}}(\bar S)$ follows from a positive scalar
curvature analogue of~\cite[Theorem F]{CaSi:emcs}.  Indeed, the proof
in~\cite[Section 6]{CaSi:emcs} generalizes to orbifolds, and reversing the
sign of the scalar curvature there, we see that any compact, connected,
totally geodesic $2$-suborbifold of an SDE orbifold of positive scalar
curvature must satisfy $\Sigma \cdot \Sigma < \chi_{\mathrm{orb}}(\Sigma)$.

\smallbreak\noindent (ii)$\Rightarrow$(iii). The positivity of the
intersection form is equivalent to the fact that $[v_j]\in \R P^1$ are in
cyclic order (see section 1), which proves the first part. Because of this,
we have $\Delta_{j-1,j}> 0$ for $1\leq j\leq k$. Now since $[\bar
S_j]\cdot[\bar S_j]<\chi_{\mathrm{orb}}(\bar S_j)$, \eqref{isn-formula}
follows from the following formulae (see section 1 for the first, the second
is standard for an orbifold $2$-sphere):
\begin{equation*}
[\bar S_j]\cdot [\bar S_j]=
\frac{\Delta_{j-1,j+1}}{\Delta_{j-1,j}\Delta_{j,j+1}},\qquad
\chi_{\mathrm{orb}} (\bar S_j)=\frac{\Delta_{j-1,j} + \Delta_{j,j+1}}
{\Delta_{j-1,j}\Delta_{j,j+1}}.
\end{equation*}
\smallbreak\noindent (iii)$\Rightarrow$(i). Under the conditions in (iii) we
are still free to cyclicly permute the $v_j$ by changing signs and
relabelling. We use this freedom to ensure that $m_0=-m_k<m_j<m_k$ for all
$0<j<k$. We want to show that conditions (a)--(b) above hold since we
already know that under these conditions $X$ admits a toric SDE metric, and
the exceptional surfaces, as fixed point sets of a subgroup of the isometry
group, are totally geodesic.

To establish (a)--(b) we note that~\eqref{isn-formula} may be rewritten as
\begin{equation}\label{isn-bis}
(n_{j+1}-n_{j})(m_j-m_{j-1})> (m_{j+1}-m_j)(n_{j}-n_{j-1}).
\end{equation}
Since $\Delta_{j-1,j}>0$, this says that $(m_{j+1},n_{j+1})$ lies on the
side containing the origin of the line $L_j$ joining $(m_{j-1},n_{j-1})$ to
$(m_j,n_j)$. We use induction to show that $(m_j)$ is increasing.  Clearly
$m_1>m_0$, so suppose $m_j>m_{j-1}$. Then dividing~\eqref{isn-bis} by
$m_j-m_{j-1}$, we see that $(m_{j+1},n_{j+1})$ is above $L_j$, hence so is
the origin. If $m_j\leq 0$, then this, together with the fact that
$\Delta_{j,j+1}>0$, shows that $m_{j+1}>m_j$. Hence the sequence
$(m_j:m_j\leq 0)$ is increasing. A similar induction starting from the fact
that $m_{k-1}<m_k$ shows that the sequence $(m_j:m_j\geq 0)$ is increasing.
Thus (a) holds, and dividing~\eqref{isn-bis} by
$(m_{j}-m_{j-1})(m_{j+1}-m_j)>0$, we obtain (b).

\smallbreak\noindent The proof of (iii)$\Rightarrow$(i) establishes the
existence of a toric SDE metric, which is unique as stated by the proof of
Theorem A and the classification of simply connected toric $4$-orbifolds.
\end{proof}
\begin{rems} The second condition in (ii) means equivalently that for any
negative toric complex structure on the complement of any fixed point in
$X$, $K_X^{-1}$ is nef. Indeed for a toric complex structure on the
complement of $x\in X$, $K_X^{-1}$ being nef is equivalent, by the
adjunction formula, to $[\bar S]\cdot[\bar S] <\chi_{\mathrm{orb}} (\bar S)$
for all $\bar S$ that do not meet $x$. Now we note that by~\cite{Joy:esd}
(see also~\cite{CaPe:emt}), $X$ admits toric scalar-flat K\"ahler metrics on
the complement of any fixed point.

Note that $X$ does not admit a global toric complex structure of either
orientation unless it is a weighted projective space (or an orbifold
quotient thereof). This can be seen by observing that toric complex
orbifolds are symplectic and so the sequence $[v_1],[v_2],\ldots [v_k]\in \R
P^1$ must have winding number two (since the $v_j$ are the normals to the
faces of a convex polytope in ${\mathfrak g}^*$). It follows easily that the
signature is $\pm (k-4)$, which equals $\mp (k-2)$ iff $k=3$.
\end{rems}

Theorem B shows that not every compact, simply connected toric $4$-orbifold
admits an SDE metric. On the other hand, as was noted in~\cite{BGMR:3s7},
there are SDE toric $4$-orbifolds with arbitrarily large second Betti
number.  It is straightforward to construct such $4$-orbifolds and compute
their rational homology using the methods presented here (the integer pairs
$(a_i,b_i)$, $i=1,\ldots k$, used in formula~\eqref{multipole} are
constrained only by pairwise linear independence).

\begin{examples} If $k=b_2(X)+2=2$ then $X$ is necessarily isometric to
an orbifold quotient of $S^4$ with the round metric: after an $\SL(2,\Z)$
transformation we may take $v_0= (-m,0)$, $v_1=(0,-n)$ and $v_2=-v_0$ as
orbifold data, $S^4$ itself being given by $mn=\pm1$.

When $k=3$, $X$ is a weighted projective space~\cite{GaLa:qro}. For
instance, orbifold data for $\C P^2$ can be taken to be $(-2,-1)$, $(-1,-1)$,
$(1,0)$, $(2,1)$.

An example with $k=4$ and only one orbifold singularity is given by the data
$(-2,-3)$, $(-1,-1)$, $(0,-1)$, $(1,0)$, $(2,3)$. More generally by taking
the $n_j$ sufficiently negative, one can construct infinitely many
examples with $b_2(X)+2=k$ for any $k\leq 2|m_0|$. The graph of $z=f_0(y)$ is
a union of line segments with integer slopes in the region $\{(y,z):z\geq
|m_0 y-n_0|\}$ as shown in Figure~\ref{fig1}.

\ifaddpics
\begin{figure}[ht]
\begin{center}
\includegraphics[width=.25\textwidth]{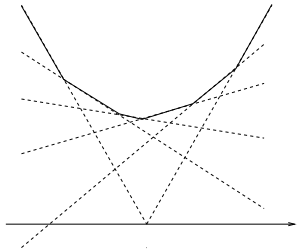}
\caption{Graph of a typical boundary value for a compact SDE orbifold.}
\label{fig1}
\end{center}
\end{figure}
\fi
\end{examples}

The body of the paper is really a series of appendices.  In section 1, we
review the classification of compact toric $4$-orbifolds,
following~\cite{OrRa:at4,HaSa:ato}. In fact, we present the classification
of simply connected compact $n$-orbifolds with a cohomogeneity two torus
action, since this is the most natural context and the fundamental paper of
Haefliger and Salem~\cite{HaSa:ato} rather understates the power of their
theory in proving results such as this.

In section 2, we review the material from~\cite{CaPe:emt,CaSi:emcs} that we
use. Then, in section 3, we present the main technical arguments that we
skipped in the proof of Theorem A.

We assume throughout that the reader is familiar with the theory of
orbifolds.

\acknowledge We thank W. Ziller for encouraging us to write this paper.  The
authors are grateful to EPSRC and the Leverhulme Trust for financial
support.  This paper was partly written while the first author was visiting
ESI and the second author was on leave at MIT. We are grateful to these
institutions for hospitality and financial support.

\numberwithin{equation}{section}
\renewcommand{\theequation}{\thesection.\arabic{equation}}

\section{Torus actions on orbifolds}

In this section we summarize the description of compact orbifolds with
torus actions due to Orlik--Raymond~\cite{OrRa:at4} and
Haefliger--Salem~\cite{HaSa:ato} (see also~\cite{HaMa:tmf,LeTo:hta}).

\subsection{Lie groups and tori acting on orbifolds}

Let $X$ be an oriented $n$-orbifold with a smooth effective action of a
compact Lie group $G$. Fix $x\in X$ with stabilizer $H\leq G$ and orbifold
structure group $\gamma$. Let $\phi\colon \tilde U \to \tilde U/\gamma= U$
be an $H$-invariant uniformizing chart about $x\in U$ and let $\tilde H$ be
the group of diffeomorphisms of $\tilde U$ which project to diffeomorphisms
induced by elements of $H$. Thus $\gamma$ is a normal subgroup of $\tilde
H$ and $H=\tilde H/\gamma$.

Elements of the Lie algebra $\mathfrak g$ of $G$ induce $\gamma$-invariant
vector fields on $\tilde U$ and the integral submanifold through
$\phi^{-1}(x)$ is $\phi^{-1}(G\cdot x\cap U)$.  Let $W=\tilde {T_x}X/\tilde
{T_x}(G\cdot x)$ be the quotient of the uniformized tangent spaces to $X$
and $G\cdot x$ at $x$. Since $\tilde H$ preserves $\tilde {T_x}(G\cdot x)$,
it acts linearly on $W$ and this induces an action of $H$ on
$W/\gamma$. Hence by the differentiable slice theorem:

\begin{quote}
{\it there is a $G$-invariant neighbourhood of the orbit $G\cdot x$ that is
$G$-equivariantly diffeomorphic to $G\times_H (B/\gamma)$, where $B$ is a
$\tilde H$-invariant ball in $W$.}
\end{quote}

Now suppose that $G=\mathfrak g/2\pi\Lambda$ is an $m$-torus ($m\leq n$),
where $\Lambda$ is a lattice in $\mathfrak g$. Then we can improve on the
above as follows. Let $U$ now be a $G$-invariant tubular neighbourhood of
$G\cdot x$ with orbifold fundamental group $\Gamma$.  Since $\pi_2(G\cdot
x)=0$, the universal orbifold cover $\pi\colon \hat U\to \hat U/\Gamma=U$
is smooth~\cite{HaSa:ato}.  Let $\hat G$ be the group of diffeomorphisms of
$\hat U$ that project to diffeomorphisms of $U$ induced by elements of $G$,
and let $\hat H$ be the stabilizer of a point $\hat x$ in $\pi^{-1}(x)$, so
that $\Gamma$ is normal in $\hat G$ and $G=\hat G/\Gamma$.  Then by the
differentiable slice theorem:
\begin{quote}
{\it there is a $G$-invariant neighbourhood of the orbit $G\cdot x$ that is
$G$-equivariantly diffeomorphic to $(\hat G\times_{\hat H} B)/\Gamma$, where
$B$ is a $\hat H$-invariant ball in $W$.}
\end{quote}
Observe that $\hat H\cap\Gamma=\gamma$, so that $(\hat G\times_{\hat H}
B)/\Gamma=(\hat G/\Gamma) \times_{\hat H/\gamma} (B/\gamma)=G\times_H
(B/\gamma)$ as before.

Since $\hat U$ is $1$-connected, $\hat G/\hat H$ is the universal cover of
$G/H$, namely $\mathfrak g/\mathfrak h$. Thus $\hat G/\hat G_0=\hat H/\hat
H_0$ is a finite group $D$, where $\hat G_0= \mathfrak g/2\pi \Lambda=G$
and $\hat H_0=\mathfrak h/2\pi\Lambda_0$ denote the identity components,
$\Lambda_0$ being a subgroup of $\Lambda$. Since $\hat H$ is the (unique)
maximal compact subgroup of $\hat G$ we have the following.
\begin{prop} \cite{HaSa:ato} Let $G=\mathfrak g/2\pi\Lambda$ be an $m$-torus
acting effectively on an oriented $n$-manifold $X$ and let $G\cdot x$ be an
orbit with $k$-dimensional stabilizer $H$. Then there is
\begin{bulletlist}
\item a rank $k$ sublattice $\Lambda_0$ of $\Lambda$,
\item a finite group $D$ with a central extension
\begin{equation*}
0\to\Lambda/\Lambda_0\to \Gamma\to D\to 1,
\end{equation*}
\item and a faithful representation $\hat H\to SO(n-m+k)$, where $\hat H$
is the maximal compact subgroup of the pushout extension $\hat
G=\Gamma\times_{\Lambda/\Lambda_0} \mathfrak g/2\pi \Lambda_0$,
\end{bulletlist}
such that a $G$-invariant tubular neighbourhood $U$ of $G\cdot x$ is
$G$-equivariantly diffeomorphic to $(\hat G\times_{\hat H} B)/\Gamma$ for a
ball $B\subset \R^{n-m+k}$. These data classify tubular neighbourhoods of
orbits up to $G$-equivariant diffeomorphism.
\end{prop}

This result is easy to apply when $k=n-m$ or $k=n-m-1$, when $\hat
H_0=\mathfrak h/2\pi\Lambda$ is a maximal torus in $SO(2(n-m))$ or
$SO(2(n-m-1)+1)$. Then $\hat H=\hat H_0$ (since $\hat H$ is in the
centralizer of $\hat H_0$), so $D=1$ and $\Gamma=\Lambda/\Lambda_0$. Hence
a tubular neighbourhood $U$ of such an orbit is classified by a subgroup
$\Lambda_0$ of $\Lambda$ such that $\hat H=\mathfrak h/2\pi\Lambda_0$.
\begin{bulletlist}
\item When $k=n-m$, $U/G$ is homeomorphic to $[0,1)^{n-m}$ and
$\Lambda_0=\bigoplus_{j=1}^{n-m}\Lambda_0^{\smash[t]{j}}$, where $
\Lambda_0^{\smash[t]{j}}$ are linearly independent rank one sublattices of
$\Lambda$ such that
$(\Lambda_0^{\smash[t]{j}}\tens_\Z\R)/2\pi\Lambda_0^{\smash[t]{j}}$ is the
stabilizer the orbits over the $j$th face of $U/G$.
\item When $k=n-m-1$, $U/G$ is homeomorphic to $[0,1)^{n-m-1}\times (-1,1)$
and $\Lambda_0=\bigoplus_{j=1}^{n-m-1}\Lambda_0^{\smash[t]{j}}$, with
$\Lambda_0^{\smash[t]{j}}$ as before.
\end{bulletlist}

To obtain a global classification, one must patch together such local
tubes. This is conveniently encoded by the {\v C}ech cohomology of $W=X/G$
with values in $\Lambda$.

\begin{prop} \cite{HaSa:ato} Suppose $W=\bigcup_i W_i$ is a union of open
sets and $(X_i,\pi_i\colon X_i\to W_i)$ are $G$-orbifolds with orbit maps
$\pi_i$.  Then there is a $G$-orbifold $(X,\pi\colon X\to W)$ with
$\pi^{-1}(W_i)$ $G$-equivariantly diffeomorphic to $X_i$ if and only if a
{\v C}ech cohomology class in $H^3(W,\Lambda)$ associated to
$\{(X_i,\pi_i)\}$ vanishes. If this is the case then the set of such
$G$-orbifolds $(X,\pi)$ is an affine space modelled on $H^2(W,\Lambda)$.
\end{prop}

\subsection{Cohomogeneity two torus actions on orbifolds}

Let us now specialize to the case $\dim W=2$ (i.e., $n=m+2$).  The union
$X_0$ of the $m$-dimensional orbits is the dense open subset on which the
action of $G$ is locally free, hence $W_0=X_0/G$ is a $2$-orbifold. The
remaining orbits have dimension $m-1$ or $m-2$, i.e., stabilizers of
dimension $k=n-m-1=1$ or $k=n-m=2$. Hence we can obtain a global
classification in this case.  (Similar arguments give a global
classification when $\dim W=1$.)

\begin{thm} \cite{OrRa:at4,HaSa:ato}
\begin{numlist}
\item Let $X$ be a compact connected oriented $(m+2)$-orbifold with a
smooth effective action of an $m$-torus $G=\mathfrak g/2\pi\Lambda$.  Then
$W=X/G$ is a compact connected oriented $2$-orbifold with boundary and
corners, equipped with a labelling $\Lambda_0^{\smash[t]{j}}$ of the edges
of $W$ \textup(the connected components of the smooth part of the
boundary\textup) by rank $1$ sublattices of $\Lambda$ such that at each
corner of $W$ the corresponding two lattices are linearly independent.

\item For any such data on $W$, there is a $G$-orbifold $X$ inducing these
data, and if $H^2(W,\Lambda)=0$ then $X$ is uniquely determined up to
$G$-equivariant diffeomorphism.

\item The orbifold fundamental group of $X$ is determined by the
long exact sequence:
\begin{equation*}
\pi_2^{orb}(W_0)\to\Lambda/{\textstyle\sum_j}\Lambda_0^{\smash[t]{j}}
\to\pi_1^{orb}(X)\to\pi_1^{orb}(W_0)\to 1.
\end{equation*}
In particular $X$ is simply connected if and only if either $W_0$ is a
smooth open disc and the lattices $\Lambda_0^{\smash[t]{j}}$ generate
$\Lambda$, or $W_0=W$ is a simply connected orbifold $2$-sphere \textup(so
that $\pi_2^{orb}(W)=\Z$\textup) and $\pi_2^{orb}(W_0)\to\Lambda$ is
an isomorphism \textup(so $m=1$\textup).
\end{numlist}
\end{thm}

\subsection{Compact toric $4$-orbifolds}

We now apply the preceding result when $m=2$ and $X$ is a simply connected
$4$-orbifold.  Then $W$ is a smooth polygonal disc with rank $1$ sublattices
$\Lambda_0^{\smash[t]{j}}\subset\Lambda\cong \Z^2$ ($j=1,\ldots k$)
labelling the edges $C_j$ of $W$, which we order cyclicly.
$\Lambda_0^{\smash[t]{j}}$ is determined by one of its generators
$v_j=(m_j,n_j)\in \Z^2$, which is unique up to a sign.  The corner
conditions mean that $v_{j-1}$, $v_j$ are linearly independent, or
equivalently $\Delta_{j-1,j}\neq 0$, for $j=1,\ldots k$ (where $v_0=-v_k$
and $\Delta_{i,j}:= m_in_j - m_jn_i$). The simple connectivity of $X$ means
that $\{v_j:j=1,\ldots k\}$ spans $\Z^2$. Since $H^2(W,\Z^2)=0$, $X$ is
uniquely determined by these data. (The classification for manifolds is more
subtle, since then we must have $\Delta_{j-1,j}=\pm 1$, i.e., adjacent pairs
of labels form a $\Z$-basis, which leads to combinatorial problems.)

We end by discussing the rational homology of such a toric $4$-orbifold
$X$.  Since $X$ is oriented and simply connected, this amounts to describing
$H_2(X,\Q)$ and its intersection form. We have already remarked that
$b_2(X)=k-2$ (and this is easy to establish by a spectral sequence
argument): in fact the closures $\bar S_j$ of the exceptional surfaces
$S_j=\pi^{-1}(C_j)$, once oriented, define rational homology classes
generating $H_2(X,\Q)$.  Obviously the only nontrivial intersections are the
self-intersections and the intersections of adjacent $\bar S_j$.  For the
latter, we note that in the orbifold uniformizing chart of order
$|\Delta_{j,j+1}|$ about $\bar S_{j}\cap \bar S_{j+1}$, the intersection
number is $\pm 1$ and hence $[\bar S_j]\cdot[\bar S_{j+1}]=\pm
1/\Delta_{j,j+1}$. Similarly, by considering the link of $S_j$ (which is an
orbifold lens space), we find that $[\bar S_j]\cdot [\bar S_j]= \pm
\Delta_{j-1,j+1}/(\Delta_{j-1,j}\Delta_{j,j+1})$.  In fact our orientation
conventions give
\begin{equation}
[\bar S_j]\cdot [\bar S_j]=
\Delta_{j-1,j+1}/(\Delta_{j-1,j}\Delta_{j,j+1}),\qquad
[\bar S_j]\cdot[\bar S_{j+1}]
= [\bar S_{j+1}]\cdot[\bar S_j]
= - 1/\Delta_{j,j+1}.
\end{equation}
We notice that $\sum_{j=1}^k m_j [\bar S_j]$ and $\sum_{j=1}^k n_j [\bar
S_j]$ have trivial intersection with any $[\bar S_i]$. Since the latter
classes span the rational homology, and the intersection form is
nondegenerate, we have $\sum_{j=1}^k m_j [\bar S_j]=0=\sum_{j=1}^k n_j [\bar
S_j]$. Since $b_2(X)=k-2$, these span the relations amongst the rational
classes $[\bar S_j]$.

For Theorem B we need a formula for the signature of $X$ (i.e., of the
intersection form on $H_2(X,\Q)$) which was given by Joyce~\cite{Joy:esd} in
the manifold case and by Hattori--Masuda~\cite{HaMa:tmf} in general. For
this formula, choose an arbitrary vector $v=(m,n)\in\Z^2$ which is not a
multiple of any $v_j$ and define $\Delta_j=m n_j-n m_j$. Then
\begin{equation*}
\sigma(X) = \sum_{j=1}^k \sign(\Delta_{j-1}\Delta_{j-1,j}  \Delta_j).
\end{equation*}
This is evidently independent of the sign choices for the $v_j$, and it is
independent of the choice of $v$, since if we move $v$ so that one
$\Delta_j$ changes sign then only two terms in the above sum change sign,
but they have the opposite sign. Let us choose the signs of the $v_j$ so
that $\Delta_j>0$ for $i=1,\ldots k$ (so that the $v_j$ lie in a half-space
bounded by the span of $v$); then $\Delta_0<0$. It follows that
$|\sigma(X)|=k-2=b_2(X)$ iff $\Delta_{j-1,j}$ all have the same sign for
$j=1,\ldots k$, i.e., iff $[v_1],\ldots [v_k]$ are in cyclic order in $\R
P^1$. Thus $\{\pm v_j:j=1,\ldots k\}$ are the normals to a
compact convex polytope in ${\mathfrak g}^*$ symmetric under $v\mapsto -v$,
as in Anguelova--Lazaroiu~\cite{AnLa:mtc}.

\section{Toric selfdual Einstein metrics}

In this section we give a brief account of the relevant results of
\cite{CaPe:emt} and \cite{CaSi:emcs}.

\subsection{Local classification}

Let $\cH^2$ denote the hyperbolic plane, which we regard as the
positive definite sheet of the hyperboloid $\{a\in S^2\R^2:\det a=1\}$
in the space $S^2\R^2$ of symmetric $2\times 2$ matrices, with the
induced metric ($-\det$ is the quadratic form of a Minkowski metric on
this vector space). 

Let $G$ denote the standard $2$-torus $\R^2/2\pi \Z^2$, with linear
coordinates $z = (z_1,z_2)$.  Consider the metric $g_F$ constructed on
(open subsets of) $\cH^2\times G$ through the formula
\begin{equation}\label{metric}
g_F = \frac{|\det\Phi|}{F^2}
\left(g_{\cH^2} +  dz\, \Phi^{-1} A \Phi^{-1} dz^t \right)
\end{equation}
where $F$ is an eigenfunction of the hyperbolic Laplacian,
\begin{equation} \label{2.22.5}
\Delta_{\cH^2} F = \frac{3}{4}F,
\end{equation}
$d z = (d z_1,d z_2)$, $A$ is the tautological $S^2\R^2$-valued
function on $\cH^2$ (with $A_a=a$), and $\Phi=\frac 12 F A-dF$ (with
$dF_a\in T^*_a\cH^2\cong a^\perp\subset S^2\R^2$).

\begin{thm} \cite{CaPe:emt} $g_F$ is an selfdual Einstein metric whose
scalar curvature scalar curvature has the opposite sign to the
quantity $\det\Phi =\frac1 4 F^2 - |dF|^2$. \textup(The metric has
singularities where $F=0$ or $\det\Phi=0$.\textup) Furthermore, any
SDE metric with nonzero scalar curvature and a $2$-torus in its
isometry group is obtained locally from this construction.
\end{thm}

A more explicit form of the metric can be obtained by introducing
half-space coordinates $(\rho(a)>0,\eta(a))$ on $\cH^2$. The standard
basis of $\R^2$ leads to a preferred choice
\begin{equation*}
A(\rho,\eta)=\frac1\rho\left[\begin{matrix}1&\eta\\
\eta&\rho^2+\eta^2\end{matrix}\right].
\end{equation*}
However---and this will be crucial later---all these local formulae are 
$SL_2(\R)$-invariant, and other half-space coordinates are given by
other unimodular bases of $\R^2$.  Identifying $\R^2$ with the Lie
algebra of $G$, it follows that we can work with any oriented
$\Z$-basis for the lattice $\Z^2$: the above formulae will transform
naturally under $SL_2(\Z)$.

One easily computes $\det A=1$ and $\det dA=(d\rho^2+d\eta^2)/\rho^2$,
so that $(\rho,\eta)$ so defined are half-space coordinates. It is
convenient to set
\begin{equation} \label{3.22.5}
f(\rho,\eta) = \sqrt{\rho}F(\rho,\eta),\quad
v_1 = (f_\rho,\eta f_\rho-\rho f_\eta),\quad
v_2 = (f_\eta,\rho f_\rho +\eta f_\eta -f),
\end{equation}
so that $\Phi=\lambda_1\tens v_1+\lambda_2\tens v_2$, where
$\lambda_1=(\sqrt\rho,0)$ and $\lambda_2 =
(\eta/\sqrt\rho,1/\sqrt\rho)$ form an orthonormal frame
($A=\lambda_1^2+\lambda_2^2$). A straightforward computation then
gives
\begin{equation}\label{1.22.5}
g_F = \frac{\rho \bigl|\ips{v_1,v_2}\bigr|}{f^2}
\left( \frac{d\rho^2 + d\eta^2}{\rho^2}
+ \frac{\ips{v_1,d z}^2 + \ips{v_2,d z}^2}{\ips{v_1,v_2}^2}
\right),
\end{equation}
where $\ips{\cdot,\cdot}$ denotes the standard symplectic form on $\R^2$.

\subsection{Boundary behaviour}

The replacement of $F$ by $f$ is quite natural because one can prove
that if $F$ satisfies \eqref{2.22.5} in a set of the form 
$\{0< \rho < a\}\times \{b<\eta < c\}$ and is of power law growth in
$\rho$, then $f$ must have an asymptotic expansion of the form
\begin{equation}\label{e1.28.4}
f(\rho,\eta)=\sqrt{\rho} F(\rho,\eta) \sim \bigl(f_0(\eta)+ f_1(\eta)\rho^2
+\cdots\bigr) +\bigl(-\tfrac12 f_0''(\eta)\rho^2+\cdots\bigr)\log\rho 
\end{equation}
where $f_0$ and $f_1$ are in general distributions.  The coefficients of the
higher powers of $\rho$ are also distributions, uniquely determined by $f_0$
and $f_1$; only even powers of $\rho$ can occur. Such boundary regularity
results are discussed in a much more general setting in \cite{Maz:hcc,
MaMe:mer}.

In addition to this local regularity, we shall need the
following uniqueness result.

\begin{prop} If $F$ satisfies \eqref{2.22.5} globally on $\cH^2$ and
$f_0=0$ on the boundary $\RP1$ of $\cH^2$, then $F=0$.
\label{van}\end{prop}

Some care is needed in interpreting this result in half-space
coordinates. At first sight, the function $F=\rho^{3/2}$ appears to be a
counter-example. However, if one changes to coordinates
\begin{equation}\label{e1.29.4.4}
\tilde\rho=\rho/(\rho^2+\eta^2), \quad\tilde\eta=-\eta/(\rho^2+\eta^2),
\end{equation}
then $F=
\tilde{\rho}^{3/2}(\tilde{\rho}^2+\tilde{\eta}^2)^{-3/2}$ and
$\sqrt{\tilde\rho}F \to 2\delta(\tilde{\eta})$ as $\tilde{\rho}\to 0$,
so that $f_0$ does not vanish at $\infty$ in the original coordinates.

More generally, if we define
$$
\tilde{f}_0(\tilde{\eta}) = \lim_{\tilde{\rho}\to 0} 
\sqrt{\tilde{\rho}}F(\rho,\eta)
$$
with $(\rho,\eta)$ and $(\tilde{\rho},\tilde{\eta})$ related as in
\eqref{e1.29.4.4}, then we see that
\begin{equation}\label{e2.29.4.4}
\tilde{f}_0(-1/\eta) = |\eta|f_0(\eta)
\end{equation}
Thus $f_0$ is the restriction of a distributional section of the
line-bundle 
$\cO(1)\tens L$ over $\RP1$, where $L$ is
the M\"obius bundle and $\cO(1)$ is the dual of the tautological line
bundle. Sections of this bundle can also be viewed as functions $\hat
f_0$ on $\R^2\punc 0$ satisfying $\hat f_0(\lambda v)= |\lambda| \hat
f_0(v)$.

An elementary way to prove Proposition~\ref{van} is via the maximum
principle.

\begin{prop} Suppose that $F$ is defined in $\cH^2$ and satisfies 
$$
\Delta F = \alpha(\alpha+1)F,\mbox{ where }\alpha \in \R.
$$
If the boundary value $f_0(\eta) = \lim_{\rho\to 0}
\rho^{\alpha}F(\rho,\eta)$ vanishes for all $\eta$ in $\RP1$, then
$F=0$.
\end{prop}
\begin{proof} We pass to the Poincar\'e model of $\cH^2$: the unit
disc with coordinates $(x,y)$, $r^2 = x^2 +y^2 < 1$. Define
$$
u = \frac{1 - r^2}{2}\mbox{ so }\Delta = 
u^2(\del_x^2 + \del_y^2).
$$
If $f = u^\alpha F$, then we have
$$
\Delta F = \Delta(u^{-\alpha}f) 
= (\Delta u^{-\alpha}) f + 2\nabla u^{-\alpha}\cdot \nabla f
+ u^{-\alpha} \Delta f.
$$
Rearranging this,
$$
\Delta f + 2u^{\alpha}
\nabla u^{-\alpha}\cdot \nabla f= \alpha(\alpha +1)f - 
u^{\alpha}(\Delta u^{-\alpha})f.
$$
By an easy computation,
$$
u^{\alpha}\Delta u^{-\alpha} = \alpha(\alpha+1) - 2\alpha^2 u.
$$
Hence
$$
\Delta f+ 2u^\alpha\nabla u^{-\alpha}\cdot \nabla f = 2\alpha^2 u f.
$$
Since $u\geq 0$ in $\cH^2$, it follows from the maximum principle that
$f$ can have neither a positive interior minimum nor a negative
interior maximum. Hence if $f\to 0$ at the boundary, then $F=0$. It is
clear that the boundary value of $f$ differs from $f_0$, defined using
half-space coordinates, by multiplication by a positive function. The
proof is complete.
\end{proof}

When $\alpha = 1/2$, we recover Propostion~\ref{van}.

Note that if we replace the eigenvalue $\alpha(\alpha+1)$ by
$\beta(\beta+1)$, then we obtain
$$
\Delta f+ 2u^{\alpha}\nabla u^{-\alpha}\cdot \nabla f = 
[(\beta-\alpha)(1+\alpha+\beta) +2\alpha^2 u] f.
$$
Since $u\to 0$ at the boundary, we cannot apply the maximum principle
if $(\beta-\alpha)(1+\alpha+\beta)<0$; in particular it is not
applicable for $\alpha>\beta=1/2$.

\subsection{The Poisson formula}

A Poisson formula reconstructs the eigenfunction $F$ from its boundary
value $f_0$; in half-space coordinates,
\begin{equation}\label{e2.28.4}
F(\rho,\eta) = \frac{1}{2}\int \frac{f_0(y)\rho^{3/2}\,d y}{\bigl(\rho^2 +
(\eta -y)^2\bigr)^{3/2}}.
\end{equation}
It is straightforward to check (e.g., by making the change of
variables $y=\eta+\rho x$) that
\begin{equation*}
f_0(\eta) = \lim_{\rho\to 0} \sqrt\rho F(\rho,\eta).
\end{equation*}
Again, despite appearances, equation~\eqref{e2.28.4} is really
$SL_2(\R)$-equivariant. Indeed, the kernel
$$
\frac{\rho}{\rho^2 + (\eta-y)^2}|d y|
$$ is $SL_2(\R)$ invariant for the diagonal action of $SL_2(\R)$ on
$\cH^2 \times \RP1$ and \eqref{e2.28.4} is the $(3/2)$-power of this
kernel applied to $f_0(y)|dy|^{-1/2}$, which is also $SL_2(\R)$
invariant.

We have seen that the map $\cP\colon f_0 \mapsto F$ given by \eqref{e2.28.4}
is injective. In fact, its image (operating on $\cD'(\RP1)$) is the space of
solutions of \eqref{2.22.5} that grow at most exponentially with geodesic
distance from a point \cite{Lew:ess}.  We shall not need this; the
interested reader is referred to Theorem~4.24 of the Introduction in
Helgason's book~\cite{Hel:gga}.

We end by remarking that~\cite[\S5.1]{CaSi:emcs} gives an integral formula
for the determinant of the matrix $\Phi=\frac 12 F A-dF$:
\begin{equation*}
\det\Phi(\rho,\eta)=-\frac{1}{8}\iint
\frac{(y-z)(\mu(y)\nu(z)-\mu(z)\nu(y))\rho^3}
{\bigl(\rho^2+(\eta-y)^2\bigr){}^{3/2}\bigl(\rho^2+(\eta-z)^2\bigr){}^{3/2}}
\,dy\,dz,
\end{equation*}
where $\mu(y)=f_0'(y)$ and $\nu(y)=yf_0'(y)-f_0(y)$.  Substituting for
$\mu$ and $\nu$, we see that
\begin{equation*}
\det\Phi(\rho,\eta)=-\frac{1}{4}\iint
\frac{(y-z)f_0(y)f_0'(z)\rho^3}
{\bigl(\rho^2+(\eta-y)^2\bigr){}^{3/2}\bigl(\rho^2+(\eta-z)^2\bigr){}^{3/2}}
\,dy\,dz.
\end{equation*}
Integrating by parts with respect to $z$ (differentiating $f_0'(z)$), we
obtain formula~\eqref{lovelyphieqn}.

\section{Boundary behaviour of $F$}

Let $X$ be a compact simply connected $4$-orbifold with an effective action
of a $2$-torus $G$, equipped with a $G$-invariant SDE metric of positive
scalar curvature. We have seen that $X/G$ is a polygonal disc and the
interior of $X/G$ is equipped with a hyperbolic metric and a hyperbolic
eigenfunction $F$. In this section we prove:
\begin{bulletlist}
\item the edges of the polygon are at infinity with respect to the
hyperbolic metric;
\item there are half-space coordinates $(\rho,\eta)$ such that if
$f_0(\eta)=\lim_{\rho\to 0}\sqrt{\rho}F(\rho,\eta)$ then on each edge
$C_j$, labelled by $\pm(m_j,n_j)$, $f_0(\eta)$ is equal, up to sign, to the
linear function $m_j\eta-n_j$;
\item $f_0(\eta)$ is continuous at the vertices of the polygon.
\end{bulletlist}
As we have remarked in section 2, for the second of these facts we can make
a unimodular change of basis and suppose that $(m_j,n_j)= (0,\ell_j)$,
where $\ell_j=\gcd(m_j,n_j)$ is the order of the orbifold structure group
of points in the corresponding special orbits. We then show that for the
half-space coordinates corresponding to such a unimodular basis,
$f_0(\eta)=\pm\ell_j$ on $C_j$.

At the corner $\bar C_j\cap \bar C_{j+1}$ corresponding to a fixed point
$x$, we would like to use the primitive vectors $(m_j,n_j)/\ell_j$ and
$(m_{j+1},n_{j+1})/\ell_{j+1}$ as a basis for $\Z^2$; unfortunately they
only form a basis for a sublattice of index $(m_j n_{j+1}-m_{j+1}n_j)
/(\ell_j\ell_{j+1})$. However, to prove the continuity of $f_0$ at the
corner, we may as well pass to the orbifold covering of a neighbourhood of
$x$ defined by this sublattice.  Hence there is no loss in supposing that
this index is $1$.

\subsection{Exceptional surfaces}\label{fermi}

Let $C$ be an edge of the polygon $W$ and $S$ its inverse image in $X$
(whose closure is an orbifold $2$-sphere). We let $2\pi/\ell$ be the cone
angle of $S$ in $X$ so that points in $S$ have orbifold structure group
$\Z_\ell$.

Near any point of $S$ we can introduce {\em Fermi coordinates}. For
$x$ near $S$, we write $r(x)$ for the distance from $x$ to $S$ and
introduce an angular coordinate $\theta$ (of period $2\pi$) such that
$d\theta$ vanishes on the radial geodesics and evaluates to $1$ on the
generator $K$ of the action of the stabilizer of $S$ (note that
$\theta$ is far from unique).  The metric then takes the form
\begin{equation}\label{fermi-metric}
g = dr^2 + r^2d\theta^2/\ell^2 + h_1 + rh_2 + r^2h_3
\end{equation}
where $h_1$ is the `first fundamental form' (restriction of the metric
to $S$), $h_2$ is the second fundamental form, and $h_3$ is a form on
$TX$ bilinear in $rd\theta$ and $TS$. Since $S$ is a fixed point set
of the isometry group generated by $K$, $h_2=0$, which we shall use in
the next subsection, but not here. (Also, by Gauss's lemma, $h_3$ does
not contain terms in $dr$, but we shall not need this precision.)

In our case, we also know that the metric is a toric SDE metric of
positive scalar curvature and so is given explicitly by
\begin{equation}\label{gF-metric}
g_F=\frac{d\rho^2+d\eta^2}{k^2 f^2}
+\frac{k^2}{f^2}
(d\psi_1,d\psi_2) P^tP
\begin{pmatrix} d\psi_1\\ d\psi_2\end{pmatrix}
\end{equation}
where
\begin{equation*}
P =\begin{pmatrix} \rho f_\eta -\eta f_\rho & f_\rho\\
f - \rho f_\rho -\eta f_\eta & f_\eta
\end{pmatrix}\qquad\text{and}\qquad k=\frac{
\sqrt{\rho}}{\sqrt{ff_\rho-\rho(f_\rho^2+f_\eta^2)}}
=\frac{\sqrt{\rho}}{\sqrt{-\det P}}.
\end{equation*}
Here the angular coordinates $(\psi_1,\psi_2)\colon X_0\to
\R^2/2\pi\Z^2$ are canonically defined up to a change of $\Z$-basis,
and the metric is invariant under such changes provided we make the
corresponding change of half-space coordinates $(\rho,\eta)$ (see
section 2). We use this freedom to let $d\psi_1$ vanish on $K$, so
that we can take $\theta=\psi_2$. There is still the freedom to add a
multiple of $\psi_1$ to $\psi_2$ and $\theta$, and we use this (in a
rather mild way) to ensure that the half-space coordinates $(\rho,\eta)$
are bounded near $S$. Note that \textit{a priori} the coordinates
$(\rho,\eta)$ are only independent in a punctured neighbourhood of $S$.

We complete the coordinates $r,\psi_1,\psi_2=\theta$ by a coordinate
$y$ so that to leading order in $r$, as a bilinear form in
$dr,dy,d\psi_1,rd\psi_2$, the metric~\eqref{fermi-metric} is given by
\begin{equation} \label{approx-metric}
dy^2 + dr^2 + a^2d\psi_1^2 + r^2d\psi_2^2/\ell^2
\end{equation}
where $a(y)>0$. The equality of the angular parts of the metrics
\eqref{gF-metric} and \eqref{approx-metric}
now reduces to the following equation:
\begin{equation*}
\frac{k^2}{f^2}P^tP =
\begin{pmatrix} a^2 & 0\\ 0 & r^2/\ell^2\end{pmatrix} +
\begin{pmatrix} O(r) & O(r^2)\\ O(r^2) & O(r^3)\end{pmatrix}.
\end{equation*}
It follows that
\begin{equation*}
\frac{k^2}{f^2} \begin{pmatrix} a^{-1} & 0\\ 0 & \ell/r \end{pmatrix}
P^tP \begin{pmatrix} a^{-1} & 0\\ 0 & \ell/r \end{pmatrix} =
I + O(r).
\end{equation*}
This matrix is symmetric and so, by binomial series expansion (for $r$
sufficiently small), it has an inverse square root which is also
symmetric and of the form $I+O(r)$.  Multiplying on both sides by this
inverse square root, we deduce that there is an orthogonal matrix with
determinant $-1$
\begin{equation*}
\begin{pmatrix} c & s\\ s & -c\end{pmatrix},\quad c^2+s^2=1
\end{equation*}
such that
\begin{equation*}
\frac{k}f P \begin{pmatrix} a^{-1} & 0\\ 0 & \ell/r \end{pmatrix}
= \begin{pmatrix} c & s\\ s & -c\end{pmatrix} +O(r),
\end{equation*}
or in other words,
\begin{equation}\label{eq-main}
\frac{k}{f} \begin{pmatrix} 
\rho f_\eta -\eta f_\rho & f_\rho\\
f- \rho f_\rho -\eta f_\eta & f_\eta
\end{pmatrix}
=  \begin{pmatrix} c a & sr/\ell \\ sa & -cr/\ell \end{pmatrix}
+ \begin{pmatrix} O(r) & O(r^2)\\ O(r) & O(r^2)\end{pmatrix}.
\end{equation}
This equation contains all the information we need.  Taking determinants, we
obtain
\begin{equation}\label{eq-det}
\frac{\rho}{f^2} = \frac{k^2}{f^2}\det P = \frac{ar}{\ell} + O(r^2).
\end{equation}
Also, if we use the $(1,2)$ and $(2,2)$ components to eliminate
$k f_\rho/f$ and $k f_\eta/f$ from the $(1,1)$ and $(2,1)$ components
of~\eqref{eq-main}, we have
\begin{align}
-\rho c r/\ell-\eta sr/\ell &= ca + O(r)\\
k - \rho sr/\ell + \eta cr/\ell &= sa + O(r).
\end{align}
In particular $c$ is $O(r)$ and so $s=\pm 1+O(r^2)$.  Thus $k\to \pm a$ as
$r\to 0$, which is nonzero.

\begin{prop} $f$ is bounded away from zero and infinity as $r\to 0$ and
$\rho$ is a defining function for the exceptional surface $S$ with
$(\rho,\eta)$ independent in a neighbourhood of $S$.
\end{prop}
\begin{proof} From the non-angular part of the metric $g_F$,
we deduce that $d\rho\wedge d\eta/f^2 = k^2 dr\wedge dy + O(r)$.
Since $d\rho\wedge d\eta$ is well defined as $r\to 0$ and $k^2
dr\wedge dy$ does not vanish at $r=0$, $f$ must be bounded. Now
by~\eqref{eq-det},
\begin{equation}\label{eq:rho-zero}
\rho=af^2\ell^{-1}r +O(r^2)
\end{equation}
and so $\rho=O(r)$. Now by the boundary regularity results discussed
in \S2, $f$ as a function of $\rho$ must be
asymptotically $O(1)$ or $O(\rho^2)$ and it follows easily
from~\eqref{eq:rho-zero} that only the first case is possible. Thus
$f$ is bounded away from zero and infinity and therefore $\rho$ has
precisely order $r$ and $d\rho\wedge d\eta$ does not vanish as
$\rho\to 0$.
\end{proof}

It is now straightforward to see that $f$ has the required boundary
behaviour.
\begin{prop} $f(\rho,\eta)= \pm \ell + O(\rho^2)$.
\end{prop}
\begin{proof} From~\eqref{eq-main}, we have $(\log f)_\rho=(\log f)_\eta
=O(\rho)$ and so $f(\rho,\eta)=A + O(\rho^2)$ for some constant
$A$. Now~\eqref{eq:rho-zero} gives $\rho=a A^2 \ell^{-1}r + O(r^2)$,
so that $d\rho=a A^2 \ell^{-1}dr+O(r)$ and $d\rho^2 = a^2A^4\ell^{-2}
dr^2+O(r)$. By comparing~\eqref{gF-metric}
and~\eqref{approx-metric} we deduce that $a^2A^4\ell^{-2}=
k^2f^2+O(r)= a^2A^2+O(r)$ so that $A^2/\ell^2=1$ (since both are
constant).
\end{proof}

We end this section by remarking that if we had interchanged the roles
of $\psi_1$ and $\psi_2$ in the above argument (corresponding to an
inversion of half-space coordinates), then we would have, in place
of~\eqref{eq-main},
\begin{equation}\label{eq-maininv}
\frac{k}{f} \begin{pmatrix} 
\rho f_\eta -\eta f_\rho & f_\rho\\
f- \rho f_\rho -\eta f_\eta & f_\eta
\end{pmatrix}
=  \begin{pmatrix} c r/\ell & sa \\ sr/\ell & -ca \end{pmatrix}
+ \begin{pmatrix} O(r^2) & O(r)\\ O(r^2) & O(r)\end{pmatrix}
\end{equation}
and hence
\begin{align}
-\rho ca -\eta sa &= cr/\ell + O(r)\\
k - \rho sa + \eta ca &= sr/\ell + O(r).
\end{align}
Assuming (as we may) that $\rho\to 0$ and $\eta$ is bounded away from
zero as $r\to 0$, the argument goes through in a similar way.  This
time, $s=O(r)$, hence $c=\pm 1+O(r^2)$ and $k\to \mp\eta a$ as $r\to
0$.  Therefore~\eqref{eq-maininv} gives $(\log f)_\rho=O(\rho)$ but
$(\log f)_\eta=1/\eta +O(\rho)$, so $f(\rho,\eta)=A \eta+O(\rho^2)$
for some constant $A$, and we compute as before that $A=\pm\ell$.

Thus $f_0(\eta)=\pm\ell\eta$, as we would expect from the coordinate
invariance of our formulae.

\subsection{Corner behaviour}

We now show that $f_0$ is continuous at a corner $\bar C_1\cap \bar C_2$,
for edges $C_1,C_2$ corresponding to exceptional surfaces $S_1,S_2$.  To do
this, we argue as in \S\ref{fermi}, but with the metric expanded about the
point $\bar S_1\cap\bar S_2$. It is now natural to introduce coordinates
$r_1$ and $r_2$, the distance functions from $S_1$ and $S_2$ respectively.
Since $S_1$ and $S_2$ are totally geodesic, the restriction of $r_2$ to
$S_1$ is also the distance function from $\bar S_1\cap \bar S_2$ and
similarly for $r_1$ and $S_2$. Therefore, after adapting the basis of
$\R^2/2\pi\Z^2$, we can suppose that the metric~\eqref{gF-metric}, to
leading order in $r_1r_2$ (as a bilinear form in $dr_1,dr_2,r_1d\psi_1,
r_2d\psi_2$), is equal to
\begin{equation*}
dr_1^2 + dr_2^2 + r_1^2d\psi_1^2/\ell_1^2 + \ell_2^2 r_2^2d\psi_2^2\ell_2^2
\end{equation*}
for some positive integers $\ell_1$ and $\ell_2$.  If we carry through
the calculations of \S\ref{fermi} with this metric, but now take into
account that $S_1$ and $S_2$ are totally geodesic (since we now need
better control over the error terms), then we certainly have
\begin{equation}
\frac{k}{f} \begin{pmatrix} 
\rho f_\eta -\eta f_\rho & f_\rho\\
f- \rho f_\rho -\eta f_\eta & f_\eta
\end{pmatrix}
=  \begin{pmatrix} c r_1/\ell_1 & sr_2/\ell_2 \\ s r_1/\ell_1 & -cr_2/\ell_2
\end{pmatrix}
+ O(r_1^2r_2^2).
\end{equation}
Hence
\begin{align}
\rho cr_2/\ell_2 +\eta sr_2/\ell_2 + cr_1/\ell_1 &= O(r_1^2r_2^2)\\
k -\rho sr_2/\ell_2 + \eta cr_2/\ell_2 - sr_1/\ell_1 &= O(r_1^2r_2^2),
\end{align}
and so away from $r_1=r_2=0$ we have
\begin{align*}
s&=\pm \frac{\rho r_2/\ell_2+r_1/\ell_1}
{\sqrt{(\rho r_2/\ell_2+r_1/\ell_1)^2+(\eta r_2/\ell_2)^2}} +O(r_1^2r_2^2)\\
c&=\mp \frac{\eta r_2/\ell_2}
{\sqrt{(\rho r_2/\ell_2+r_1/\ell_1)^2+(\eta r_2/\ell_2)^2}} +O(r_1^2r_2^2)\\
k&=\pm \sqrt{(\rho r_2/\ell_2+r_1/\ell_1)^2+(\eta r_2/\ell_2)^2}+O(r_1^2r_2^2).
\end{align*}
(We must be careful as $s$ and $c$ are not continuous at $r_1=r_2=0$.)
We deduce from these formulae that
\begin{equation}\label{ub}\begin{split}
(\log f)_\rho & =\frac{\rho r_2^2/\ell_2^2+r_1r_2/(\ell_1\ell_2)}
{(\rho r_2/\ell_2+r_1/\ell_1)^2+(\eta r_2/\ell_2)^2} +O(r_1^2r_2^2)\\
(\log f)_\eta & =\frac{\eta r_2^2/\ell_2^2}
{(\rho r_2/\ell_2+r_1/\ell_1)^2+(\eta r_2/\ell_2)^2} +O(r_1^2r_2^2).
\end{split}\end{equation}
Let the corner correspond to $\eta=\eta_0$. We shall investigate the
behaviour of $f$ in polar coordinates centred at $(0,\eta_0)$, by
introducing
\begin{equation*}
\eta -\eta_0 = R\cos\Theta,\quad \rho = R\sin\Theta.
\end{equation*}
Then 
\begin{equation*}
\del_R\log f  = ({\cos\Theta}\, \del_\eta +{\sin\Theta}\,\del_\rho)\log f
\end{equation*}
and this is uniformly bounded for any fixed $\Theta$ by \eqref{ub}. Hence
$f(R,\Theta)$ has a limit $f(0,\Theta)$ as $R\to 0$, for each fixed
$\Theta\in(0,\pi)$. It also has a limit as $\Theta\to 0$ or $\pi$ for each
fixed $R>0$, namely the boundary value $f_0(\eta_0\pm R)$.  Thus $f$ is
bounded in $[0,\eps]\times[-\pi,\pi]$ for some $\eps>0$.

Now we note that
\begin{equation*}
\del_\Theta f = f \del_\Theta \log f=
f R({-\sin\Theta}\,\del_\eta +{\cos\Theta}\,\del_\rho)\log f
\end{equation*}
and the right hand side is $O(R)$, uniformly in $\Theta$ by
\eqref{ub}.  Hence by integration, we have
$|f(R,\Theta_1)-f(R,\Theta_2)|=O(R)$ for any $\Theta_1,\Theta_2\in
(0,\pi)$.  We deduce that
$$|f_0(\eta_0+R)-f_0(\eta_0-R)|=O(R)$$ and hence finally, taking $R\to
0$, that $f_0$ is continuous at the corner.

%
%
\newcommand{\bauth}[1]{\mbox{#1}} \newcommand{\bart}[1]{\textit{#1}}
\newcommand{\bjourn}[4]{#1\ifx{}{#2}\else{ \textbf{#2}}\fi{ (#4)}}
\newcommand{\bbook}[1]{\textsl{#1}}
\newcommand{\bseries}[2]{#1\ifx{}{#2}\else{ \textbf{#2}}\fi}
\newcommand{\bpp}[1]{#1} \newcommand{\bdate}[1]{ (#1)} \def\band/{and}
\newif\ifbibtex
\ifbibtex
\bibliographystyle{genbib}
\bibliography{papers}
\else

\fi

\end{document}